\documentclass[a4paper]{article}
\usepackage{amsthm,amsmath,amssymb}
\newtheorem{teor}{Theorem}[section]
\newtheorem{defin}[teor]{Definition}
\newtheorem{lemm}[teor]{Lemma}
\newtheorem{osse}[teor]{Remark}
\newtheorem{prop}[teor]{Proposition}
\newtheorem{defi}[teor]{Definition}
\newtheorem{coro}[teor]{Corollary}
\newtheorem{prob}[teor]{Problem}

\newcommand{\bele}{\begin{lemm}\begin{sl}}
\newcommand{\enle}{\end{sl}\end{lemm}}
\newcommand{\bedef}{\begin{defi}\begin{sl}}
\newcommand{\eddef}{\end{sl}\end{defi}}
\newcommand{\bete}{\begin{teor}\begin{sl}}
\newcommand{\ente}{\end{sl}\end{teor}}
\newcommand{\beos}{\begin{osse}\begin{rm}}
\newcommand{\eddos}{\end{rm}\end{osse}}
\newcommand{\bepr}{\begin{prop}\begin{sl}}
\newcommand{\empr}{\end{sl}\end{prop}}
\newcommand{\bepro}{\begin{prob}\begin{rm}}
\newcommand{\empro}{\end{rm}\end{prob}}
\newcommand{\bede}{\begin{defin}\begin{sl}}
\newcommand{\edde}{\end{sl}\end{defin}}
\newcommand{\beco}{\begin{coro}\begin{sl}}
\newcommand{\enco}{\end{sl}\end{coro}}

\newcommand{\quext}{\quad\text}
\newcommand{\qquext}{\qquad\text}
\newcommand{\de}{\partial}

\newcommand{\RR}{\mathbb{R}}
\newcommand{\NN}{\mathbb{N}}

\newcommand{\Dzero}{\mathbb{D}_0}

\newcommand{\beeq}[1]{\begin{equation}\label{#1}}
\newcommand{\eddeq}{\end{equation}}

\newcommand{\beeqa}[1]{\begin{eqnarray}\label{#1}}
\newcommand{\eddeqa}{\end{eqnarray}}

\newcommand{\beal}[1]{\begin{align}\label{#1}}
\newcommand{\eddal}{\end{align}}

\newcommand{\bespl}[1]{\begin{split}\label{#1}}
\newcommand{\edspl}{\end{split}}

\newcommand{\bega}[1]{\begin{gather}\label{#1}}
\newcommand{\edga}{\end{gather}}

\newcommand{\beeqax}{\begin{eqnarray*}}
\newcommand{\eddeqax}{\end{eqnarray*}}

\def\qed{\ifmmode 
  \else \leavevmode\unskip\penalty9999 \hbox{}\nobreak\hfill
  \fi
  \quad\hbox{\hskip.5em\vrule width.4em height.6em depth.05em\hskip.1em}}
\def\endproofsym{\qed}
\newcommand{\dimbox}{\hbox{\hskip.5em\vrule width.4em height.6em depth.05em\hskip.1em}}
\renewenvironment{proof}[1][Proof]{\trivlist\item[\hskip\labelsep{\hskip0pt
    {\normalfont\scshape#1.}\hskip .321429\parindent}]\ignorespaces}
{\endproofsym\endtrivlist}
\def\endnobox{\def\endproofsym{}\end{proof}\def\endproofsym{\qed}}

\newcommand{\no}{\nonumber}

\newcommand{\beeqao}{\begin{eqnarray}\no}
\newcommand{\bealo}{\begin{align}\no}
\newcommand{\besplo}{\begin{split}\no}
\newcommand{\begao}{\begin{gather}\no}

\newcommand{\duav}[1]{\langle{#1}\rangle}

\newcommand{\cc}{{\mathfrak c}}

\newcommand{\perogni}{\forall\,}
\newcommand{\esiste}{\exists\,}
\newcommand{\itt}{\int_0^t}

\newcommand{\io}{\int_\Omega}

\newcommand{\OO}{_{\Omega}}

\newcommand{\bn}{\boldsymbol{n}}

\newcommand{\dn}{\partial_{\bn}}

\newcommand{\lhs}{left hand side}
\newcommand{\rhs}{right hand side}

\DeclareMathOperator{\dive}{div}
\DeclareMathOperator{\deriv}{d}

\DeclareMathOperator{\dist}{dist}

\DeclareMathOperator{\discr}{discr}

\newcommand{\HUH}{H^1(0,T;H)}

\newcommand{\CZV}{C^0([0,T];V)}

\newcommand{\LDHD}{L^2(0,T;H^2(\Omega))}

\let\TeXchi\chi
\def\chi{{\setbox0 \hbox{\mathsurround0pt
$\TeXchi$}\hbox{\raise\dp0 \copy0 }}}

\newcommand{\zzn}{_{0,n}}

\newcommand{\teta}{\vartheta}

\newcommand{\calH}{{\mathcal H}}

\newcommand{\calX}{{\mathcal X}}

\newcommand{\calUl}{{\mathcal W}_\ell}
\newcommand{\calA}{{\mathcal A}}

\newcommand{\calE}{{\mathcal E}}
\newcommand{\calEq}{\mathbb{E}_0}

\newcommand{\calM}{{\mathcal M}}

\newcommand{\calB}{{\mathcal B}}

\newcommand{\dX}[4]{d_{{\mathcal X}}\big((#1,#2),(#3,#4)\big)}

\newcommand{\calV}{{\mathcal V}}
\newcommand{\calY}{{\mathcal Y}}
\newcommand{\calZ}{{\mathcal Z}}
\newcommand{\calW}{{\mathcal W}}

\newcommand{\dit}{\deriv\!t}
\newcommand{\dis}{\deriv\!s}
\newcommand{\ditau}{\deriv\!\tau}

\newcommand{\ddt}{\frac{\deriv\!{}}{\dit}}

\newcommand{\calL}{{\cal L}}

\numberwithin{equation}{section}

\begin{document}

\title{Global and 
  exponential attractors for the Penrose-Fife system}
%
%
%

%
\author{Giulio Schimperna\\
Dipartimento di Matematica, Universit\`a di Pavia,\\
Via Ferrata, 1, I-27100 Pavia, Italy,\\
E-mail: {\tt giusch04@unipv.it}
}

\maketitle
\begin{abstract}
 The Penrose-Fife system for phase transitions is addressed.
 Dirichlet boundary conditions for the temperature are
 assumed. Existence of global and exponential attractors is proved. 
 Differently from preceding contributions, here
 the energy balance equation is both singular at 0
 and degenerate at $\infty$. For this reason, the 
 dissipativity of the associated dynamical process
 is not trivial and has to be proved rather carefully.
\end{abstract}
{\bf AMS (MOS) subject clas\-si\-fi\-ca\-tion:}~~35B41, 
 35K55, 80A22

%
%
%



%
%


\section{Introduction}
\label{secintro}

We consider here the thermodynamically consistent model 
for phase transitions proposed by Penrose and Fife in
\cite{PF,PF2} and represented by the equations
\beal{calorein}
   & \teta_t+\lambda(\chi)_t 
     +\dive\Big(m(\teta)\nabla \frac1\teta\Big)=g,\\
  \label{phasein}
   & \chi_t-\Delta\chi+W'(\chi)=
    \lambda'(\chi)\Big(-\frac1\teta+\frac1\teta_c\Big).
\end{align}
The system above is settled in a smooth, bounded domain 
$\Omega\subset\RR^3$, with boundary $\Gamma$.
The unknowns are the {\sl absolute}\/ temperature $\teta>0$ 
and the order parameter $\chi$. The smooth functions $\lambda'$,
$m$ and $W$ represent the latent heat, the thermal 
conductivity, and the potential associated to the local
phase configuration, respectively, and 
$\teta_c>0$ is a critical temperature. 
Finally, $g$ is a volumic heat source. On the basis
of physical considerations, the {\sl kinetic}\/ equation
\eqref{phasein} is complemented, as usual,
with no-flux (i.e., homogeneous Neumann) boundary conditions;
instead, various types of meaningful boundary conditions 
can be associated with the {\sl energy balance}\/
equation~\eqref{calorein}. We shall consider here
the Dirichlet boundary conditions.

As far as well-posedness is concerned, system 
\eqref{calorein}--\eqref{phasein} has been studied
in a number of recent works, among which we quote
\cite{CGRS,CL,CP,GM,La1,La2,SpZ},
under various assumptions on the data.
The papers listed above also contain a
much more comprehensive bibliography.
Just a rapid survey of the literature
suggests that, indeed, the choice of the boundary 
conditions for $\teta$ can give rise to several
different mathematical situations. 
In particular, the Dirichlet and Robin 
conditions seem easier to treat than the Neumann ones
(cf., e.g., \cite{CGRS,GM} for further comments),
due to correspondingly higher coercivity.
Another important factor is the expression 
of the thermal conductivity $m$.
Meaningful choices are
given by (cf.~\cite{CL} for further comments)
\beeq{choicesm}
  m(r)\sim m_0 r + m_\infty r^2,
   \qquad m_0,m_\infty\in[0,\infty).
\end{equation}
In particular, $m_0=0,m_\infty>0$ represents the 
Fourier heat conduction
law, which appears to be the most difficult
situation \cite{La2} since equation 
\eqref{calorein}, which is now linear in $\teta$,
is coupled with the {\sl singular}\/ relation
\eqref{phasein}. Instead, in the case
$m_0>0,m_\infty=0$, the well-posedness issue
is simpler (cf.~\cite{La1,SpZ}); however, there is
a lack of coercivity for large $\teta$,
which creates difficulties in the long-time
analysis. Finally, the probably simplest situation
is that proposed in \cite{CL} (see also \cite{CLS}),
i.e., $m_0,m_\infty>0$, 
since \eqref{calorein} maintains
both the singular character at $0$ 
and the coercivity at $\infty$.

In view of these considerations, it is not 
surprising that the long time behavior of 
\eqref{calorein}--\eqref{phasein} is 
better understood when
$m_0,m_\infty>0$, and in this case
the existence of the global attractor
has been shown in \cite{RS,RS4}.
Indeed, testing \eqref{calorein} by $\teta$ 
one readily gets a dissipative estimate 
for the temperature, which permits to construct
a uniformly absorbing set and, consequently,
the global attractor. Similar results
are also obtained in \cite{IK,KK}, where
it is actually taken $m_\infty=0$, but 
a term $\mu_\infty\teta$, with $\mu_\infty>0$, 
is added on the \lhs\ of \eqref{calorein},
so that the system is still coercive 
in~$\teta$.

Speaking of the non-coercive case $m_\infty=0$,
up to our knowledge the only papers devoted
to the large-times analysis of it are
\cite{SZ} (see also \cite{SZ2} for the {\sl conserved}\/ 
case) and \cite{FS}. In \cite{SZ}, the case of homogeneous
Neumann conditions for both unknowns
is addressed in one space dimension, and 
existence of a global attractor is shown in 
a proper phase space which takes into account
the conservation (or dissipation) properties coming from
the no-flux conditions. In \cite{FS}, 
the (non-homogeneous)
Dirichlet case is considered in three
space dimensions and $\omega$-limits of 
single trajectories are studied. It is 
worth remarking
that in both papers the external
source $g$ is taken equal to $0$.

In the present work, we provide a further
contribution to the analysis of the noncoercive
case. Precisely, we assume $m_0>0$, $m_\infty=0$,
and take Dirichlet boundary conditions for $\teta$
exactly as in \cite{FS}. For the resulting problem, we show
existence of both global and exponential attractors.
Comparing with \cite{FS},
where the behavior of a single trajectory
is investigated, here the proofs are 
very different and in several
points more difficult. Indeed, determining attractors
means to understand the behavior of {\sl bundles}\/ of 
trajectories, so that we need to find estimates
which are uniform not only in time, but also
with respect to initial data varying in a bounded set.
We then try to minimize technicalities
by making some restrictions on data.
Namely, we take a constant latent heat (i.e., set
$\lambda(\chi)=\chi$), set $m(r)=1$ (i.e.,
$m_0=1$, $m_\infty=0$), let the critical
temperature $\teta_c$ be equal to 1, and correspondingly
assume the Dirichlet condition $\teta\equiv 1$
on the boundary. Actually, all these assumptions
could by avoided by paying the price of 
some additional computations in the proofs.
More restrictive is, instead, the assumption
$g=0$, which we take exactly as 
it was done in \cite{FS,SZ}.
We then end up with the system
\beal{calorein2}
   & \teta_t+\chi_t 
     -\Delta\Big(-\frac1\teta\Big)=0,\\
  \label{phasein2}
   & \chi_t-\Delta\chi+W'(\chi)=
     1-\frac1\teta.
\end{align}
Being $g=0$, \eqref{calorein2}--\eqref{phasein2}
admits a Liapounov functional (and consequently 
a dissipation integral),
and this information will be crucial 
to overcome the lack of coercivity in $\teta$.
Actually, the global attractor will be
constructed by proving uniform boundedness
and asymptotic compactness of single trajectories
and taking advantage of the dissipation property.
Although this procedure might seem 
straighforward, the proof presents
a number of difficulties. First of all,
we have to settle the problem in
a phase space $\calX$ (cf.~\eqref{defiX}
below) where both $\teta$ and $\chi$
are bounded in sufficiently strong norms. 
The conditions we require on the initial data
are in fact more restrictive 
than what is necessary, e.g., for the
mere well-posedness. In particular,
we cannot deal with completely general
potentials $W$. Namely, we are forced
to assume $W$ be a smooth function defined
on the whole real line (like, e.g., 
the {\sl double well}\/ potential 
$W(r)=(r^2-1)^2$), and, for instance,
we cannot treat the {\sl singular}\/ potentials,
i.e., those being identically $+\infty$ 
outside a bounded interval, like 
the so-called {\sl logarithmic}\/ 
potential $W(r)=(r+1)\log(r+1)+(1-r)\log(1-r)-\lambda r^2$,
where $\lambda>0$. Moreover, we note that, in
analogy with the coercive case $m_\infty>0$
studied in~\cite{RS}, $\calX$ does not have 
a Banach structure, due to the nonlinear terms in the 
energy, but it is just a metric space.
In this setting, the key point of our argument
is the proof of a uniform time regularization
property for the solutions, which, in our
opinion, can constitute an interesting issue
by itself. Namely,
we can show that both $\teta$ and $u=\teta^{-1}$ 
are uniformly bounded for sufficiently
large times, whereas this need not hold for 
the initial temperature $\teta_0$. 
Thus, \eqref{calorein2} eventually loses 
both the singular and the degenerate character.

A further open problem to which we give a positive answer
is the existence of exponential attractors for the system
\eqref{calorein2}--\eqref{phasein2}. This is shown by using
the so-called method of $\ell$-trajectories 
(cf.~\cite{MP,Pr1,Pr2,Pr3}). However, we cannot
prove exponential attraction in the metric of 
$\calX$ (that keeps, in some way, a trace 
of the nonlinear terms),
but are forced to work with a weaker norm, 
corresponding in fact to the only contractive 
estimate which seems to hold for 
system~\eqref{calorein2}--\eqref{phasein2}. 

The rest of this paper is organized as
follows. In the next Section~\ref{secmain}, we present
our hypotheses and state our main results.
The proofs are collected in Section~\ref{secproofs}.

\vspace{2mm}

\noindent%
{\bf Acknowledgment.}~~%
We express our gratitude to Elisabetta Rocca and Riccarda
Rossi for fruitful discussions on the subject of this work.


\section{Notation and main results}
\label{secmain}

Let $\Omega\subset\RR^3$ be a smooth bounded domain
with boundary $\Gamma$.
Let us set $H:=L^{2}(\Omega)$ and denote by
$(\cdot,\cdot)$ both the scalar product in $H$ and that
in $H\times H$, and by $\|\cdot\|$ the induced norm. 
The symbol $\|\cdot\|_{X}$ indicates the norm 
in the generic Banach space $X$. 
Next, we set $V:=H^1(\Omega)$, $V_0:=H^1_0(\Omega)$,
and define
\beal{defiA}
  & A:V_0\to V_0', \quad
   \duav{Av,z}_0:=\io \nabla v\cdot \nabla z, 
    \quad\perogni v,z\in V_0,\\
 \label{defiB}
  & B:V\to V', \quad
   \duav{Bv,z}:=\io \big(vz+\nabla v\cdot \nabla z\big),
    \quad\perogni v,z\in V,
\end{align}
$\duav{\cdot,\cdot}_0$ and $\duav{\cdot,\cdot}$ denoting
the duality pairings between $V_0$ and 
$V_0'=H^{-1}(\Omega)$ and 
between $V$ and $V'$, respectively. It turns out
that $A$ and $B$ are the Riesz operators associated
to the standard norms in $V_0$ and $V$, respectively.

%
Our hypotheses on the potential $W$ are the following:
%
\beal{W1}
  & W\in C^2(\RR;\RR), \quad W'(0)=0, \quad
   \lim_{r\to\infty}W'(r)r=+\infty,\\
 \label{W2}
  & \esiste \lambda\ge0:~~W''(r)\ge-\lambda~~\perogni r\in \RR.
\end{align}
In particular, by the latter assumption,
$\beta(r):=W'(r)+\lambda r$ is increasingly
monotone. 
Next, considering $B$, with a small 
abuse of notation, as a strictly positive
unbounded linear operator on $H$ with 
domain $D(B)=\{v\in H^2(\Omega):
\dn v=0\text{ on $\de\Omega$}\}$, 
we can take real powers of $B$ 
and set $\calV_{2s}:=D(B^s)$, endowed with
the graph norm $\|v\|_s:=\|B^s v\|$. Note
that $\calV_1=V$.
The variational formulation
of system~\eqref{calorein2}--\eqref{phasein2}
takes then the form
\beal{calore}
   & \teta_t+\chi_t 
     +A\Big(1-\frac1\teta\Big)=0,
   \qquext{in }\, V_0',\\
  \label{phase}
   & \chi_t+B\chi+W'(\chi)=
     1-\frac1\teta,
   \qquext{in }\, V'
\end{align}
(in order to get the Riesz map $B$, 
$\chi$ has been added and subtracted from
the \lhs, and $-\chi$ has been included in 
$W'$). Next, we define the 
associated {\sl energy}\/ functional as:
\beeq{defiE}
  \calE=\calE(\teta,\chi)
   :=\io \Big(\teta-\log\teta
    +\frac12|\chi|^2
    +\frac12|\nabla \chi|^2
    +W(\chi)\Big).
\end{equation}
We immediately observe that $\calE$ is finite and 
bounded from below on the ``energy space''
\beeq{defiXE}
  \calX_\calE:=\big\{(\teta,\chi):\teta\in L^1(\Omega),~
   \teta>0\text{ a.e.~in }\Omega,~\log\teta\in L^1(\Omega),~
   \chi\in V,~W(\chi)\in L^1(\Omega)\big\}.
\end{equation}
Nevertheless, due to the lack of coercivity 
(and consequently of compactness) in $\teta$
(the finiteness of energy only implies 
that $\teta\in L^1(\Omega)$),
no existence result is known, up to our knowledge,
for data lying just in 
$\calX_\calE$. Namely, noting as 
Problem~(P) the coupling of 
\eqref{calore}--\eqref{phase}
(intended to hold for a.e.~value of
time in $(0,\infty)$) with the 
initial condition
\beeq{iniz}
  \teta|_{t=0}=\teta_0, \quad
   \chi|_{t=0}=\chi_0,
   \qquext{a.e.~in }\,\Omega,
\end{equation}
we have the following result,
proved in \cite[Thm.~2.1]{FS} (see also 
\cite[Prop.~2.1]{GM}):
\bete\label{teoFS}
 Let \eqref{W1}--\eqref{W2} hold and let 
 $(\teta_0,\chi_0)\in \calX_\calE$. Let, in addition
 $\teta_0\in L^p(\Omega)$ for some
 $p>6/5$. Then, there exists one and only
 one couple $(\teta,\chi)$ solving\/
 {\rm Problem~(P)} and such that
 \beal{regtetaFS}
   & \teta\in H^1(0,T;H^{-1}(\Omega))\cap
    L^\infty(0,T;L^p(\Omega)), \qquad 
    \teta>0\quext{a.e.~in }\,\Omega\times(0,T),\\
  \label{reg1tetaFS}
   & \big(1-1/\teta\big)\in L^2(0,T;V_0),\\
  \label{regchiFS}
   & \chi\in\HUH\cap\CZV\cap\LDHD,
 \end{align}
 hold for all $T>0$. Such a couple will be called a\/
 {\rm ``solution''} in the sequel.
\ente
Since we need to control uniformly in time
the ``large values'' of the temperature, we have 
to ask a bit more summability on $\teta_0$
and a bit more regularity on $\chi_0$. Correspondingly,
we will also get some more regularity
than \eqref{regtetaFS}--\eqref{regchiFS}.
Namely, we set
\beeq{defiX}
  \calX:=\big\{(\teta,\chi)\in \calX_\calE:\teta\in L^p(\Omega),~
     \chi\in\calV_{\frac{3+\epsilon}2}\big\},
\end{equation}
where we assume that
\beeq{epsip}
   \epsilon\in(0,1), \qquad p>3.
\end{equation}
Actually, we need $\epsilon>0$ 
in order to ensure that $\chi$ stays
in $L^\infty(\Omega)$, while the higher summability
of $\teta_0$ seems necessary to get a uniform
in time estimate for $\teta(t)$.

We remark that the set $\calX$, which of course has
no linear structure, can be endowed with a complete metric
which makes it a suitable phase space for
the associated dynamical process.
As in \cite{RS} (see also \cite{RSS,Se}),
we can take
\bealo
  & \dX{\teta_1}{\chi_1}{\teta_2}{\chi_2}
   :=\|\teta_1-\teta_2\|_{L^p(\Omega)}
   +\|\chi_1-\chi_2\|_{\frac{3+\epsilon}2}\\
 \label{defidX}
  & \mbox{}~~~~~~~~~~
   +\|\log^-\teta_1-\log^-\teta_2\|_{L^1(\Omega)}
   +\|\beta(\chi_1)-\beta(\chi_2)\|_{L^1(\Omega)},
\end{align}
where $(\cdot)^-$ denotes negative part
(notice, however, that the latter term could
be omitted since it is dominated by the second one
due to \eqref{W1} and the continuous 
embedding $\calV_{\frac{3+\epsilon}2}\subset L^\infty(\Omega)$).
Correspondingly, we take initial data such 
that
\beeq{hpiniz}
  (\teta_0,\chi_0)\in\calX.
\end{equation}
In the sequel, we will denote by $S(t)$ the 
semigroup operator associating to $(\teta_0,\chi_0)$ 
the corresponding solution evaluated at time $t$. 
The proof that $S(\cdot)$ fulfills the usual
properties of a continuous semigroup on $\calX$
is more or less standard and can be carried out
along the lines, e.g., of~\cite[Sec.~4]{RS}. Hence, 
we omit the details. Instead, we focus on regularization
properties of $S(t)$. The key step of our 
investigation is the following
\bete\label{teounif}
 Let\/ \eqref{W1}--\eqref{W2}, \eqref{epsip}
 hold and let
 $B$ be a set of initial data bounded in 
 $\calX$. More precisely, let
 $\Dzero$ stand for the $d_\calX$-radius of the set, namely
 \beeq{Dzero}
   \Dzero:=\sup_{(\teta_0,\TeXchi_0)\in B}\dX{\teta_0}{\chi_0}10.
 \end{equation}
 Then, letting $(\teta_0,\chi_0)\in B$ and 
 $(\teta(t),\chi(t)):=S(t)(\teta_0,\chi_0)$, 
 there exist a time $T_\infty>0$ 
 and a constant $Q_\infty$ 
 depending only on $\Dzero$, such that,
 for all $t\ge T_\infty$, there holds
 \beal{rego1}
   & \|\teta(t)\|_{V\cap L^\infty(\Omega)}
    +\|\teta^{-1}(t)\|_{V\cap L^\infty(\Omega)}\le Q_\infty,\\
   \label{rego2}
   & \|\chi(t)\|_{H^2(\Omega)}\le Q_\infty.
 \end{align}
\ente
\beos\label{nuovarego}
 Suitably modifying the proofs,
 one could show that {\sl any}\/ 
 strictly positive time 
 could be taken as $T_\infty$.
 We omit the proof
 of this fact since it would involve
 further technical complications.
 We just notice that
 the quantity $Q_\infty$ in 
 \eqref{rego1}--\eqref{rego2} 
 would then depend on $T_\infty$ and explode
 as $T_\infty\searrow 0$. 
\eddos
Notice that the 
bounds~\eqref{rego1}--\eqref{rego2}
are somehow weaker than a true dissipative 
estimate. Nevertheless, they will suffice for
the proof of our main result (for the definition
of the global attractor we refer 
to the monograpgh~\cite{Te}):
\bete\label{teoattra}
 Let the assumptions
 of\/ {\rm Theorem~\ref{teounif}} hold.
 Then, the semigroup $S(\cdot)$ associated with
 {\rm Problem (P)} admits the global attractor
 $\calA$, which is compact in $\calX$. More
 precisely, 
 \beeq{compA}
   \esiste c_{\calA}>0:~~
    \|\teta\|_{V\cap L^\infty(\Omega)}
    +\|\teta^{-1}\|_{V\cap L^\infty(\Omega)}
    +\|\chi\|_{H^2(\Omega)}
   \le c_{\calA}
    \quad\perogni(\teta,\chi)\in \calA.
 \end{equation}
\ente
Finally, we can prove existence of an exponential
attractor:
\bete\label{teoespo}
 Let the assumptions
 of~{\rm Theorem~\ref{teounif}} hold.
 Then, the semigroup $S(\cdot)$ associated with
 {\rm Problem~(P)} admits an exponential attractor 
 $\calM$. More precisely, $\calM$ is a compact
 set of $\calX$, which has finite fractal 
 dimension in $V_0'\times H$, such that for any
 bounded set $B\subset \calX$ there holds
 \beeq{espoattra}
    \dist(S(t)B,\calM)\le Q(\Dzero) e^{-\kappa t},
      \qquad\perogni t\ge 0,
 \end{equation}
 where $\dist$ represents the unilateral Hausdorff
 distance of sets with respect to the (product)
 norm in $V_0'\times H$, $\kappa>0$ is independent
 of $B$, $Q$ is a monotone function, and 
 $\Dzero$ is the $\calX$-radius of $B$ given
 by~\eqref{Dzero}.
\ente
\beos\label{primaregoespo}
 As noted in the Introduction, we use
 the (rather weak) topology of $V_0' \times H$ 
 since it seems difficult to prove a contractive 
 estimate in a better norm. Further comments
 will be given at the end of the proof
 (cf.~Remark~\ref{secondaregoespo} below).
\eddos


\section{Proofs}
\label{secproofs}

In what follows, the symbols $c$, 
$\kappa$, and $c_i$, 
$i\ge0$, will denote positive constants
depending on $W,\Omega$, and independent
of the initial datum and of time. 
The values of $c$ and $\kappa$
are allowed to vary even within the same line.
Moreover, $Q:\RR^+\to\RR^+$ denotes a 
generic monotone function.
Capital letters like $C$ or $C_i$ will be used to 
indicate constant which have other dependencies
(in most cases, on the initial datum). 
Finally,
the symbol $c\OO$ will denote some embedding constants
depending only on the set $\Omega$. 

\vspace{2mm}

\noindent%
{\bf Proof of Theorem~\ref{teounif}.}~~%
The basic idea to prove the uniform
bounds \eqref{rego1}--\eqref{rego2}
is to combine an estimate in a small 
interval $[0,T_0]$, where $T_0$ depends on $\Dzero$, 
with a further uniform estimate
holding on $[T_0,\infty)$. This procedure requires
a number of steps, which are carried out below.
Notice that some parts of the procedure might
have a formal character in the present regularity
setting (e.g., test functions could be not
regular enough). However, all the procedure
could be standardly made rigorous by working
on some approximation and then passing to 
the limit (notice that the solution
is known to be unique). We omit the details
of this straighforward argument, for brevity.

\vspace{2mm}

\noindent%
{\bf First estimate.}~~%
We start by deriving the energy estimate. Testing
\eqref{calore} by $1-1/\teta$, we have
\beeq{energy1}
  \ddt\io \big(\teta-\log\teta\big)
   +\Big\|1-\frac1\teta\Big\|^2_{V_0}
   =-\Big\langle\chi_t,1-\frac1\teta\Big\rangle.
\end{equation}
Next, multiplying \eqref{phase} by $\chi_t$, we obtain
\beeq{energy2}
  \ddt\io \Big(\frac12|\chi|^2+\frac12|\nabla\chi|^2
   +W(\chi)\Big)
   +\|\chi_t\|^2
   =\Big\langle\chi_t,1-\frac1\teta\Big\rangle,
\end{equation}
whence, summing \eqref{energy1} and \eqref{energy2}
and recalling \eqref{defiE}, 
\beeq{energy}
  \ddt\calE
   +\Big\|1-\frac1\teta\Big\|^2_{V_0}
   +\|\chi_t\|^2
   =0.
\end{equation}
In particular, integrating from $0$ to an arbitrary $t>0$ we
have
\beeq{energycons}
  \calE(t)
   +\itt\Big(\Big\|1-\frac1\teta\Big\|^2_{V_0}
   +\|\chi_t\|^2\Big)
  \le \calE(0)
  \le Q(\Dzero).
\end{equation}

\vspace{2mm}

\noindent%
{\bf Second estimate.}~~%
We test \eqref{phase} by $2B^{\frac{1+\epsilon}2}\chi_t$. 
Thanks to $\epsilon<1$ and using
Poincar\'e's and Young's inequalities, we obtain
\beeq{conto21}
  \ddt\|\chi\|_{\frac{3+\epsilon}2}^2
  +\|\chi_t\|_{\frac{1+\epsilon}2}^2
  \le c\Big\|1-\frac1\teta-W'(\chi)\Big\|^2_{\frac{1+\epsilon}2}
  \le c_1\Big\|1-\frac1\teta\Big\|^2_{V_0}
   + c \|W'(\chi)\|_{V}^2.
\end{equation}
Then, we note that,
by \eqref{W1} and the continuous embedding 
$\calV_{\frac{3+\epsilon}2}\subset L^\infty(\Omega)$, 
\beeq{conto22}
  \|W'(\chi)\|_{V}^2
   = \|W'(\chi)\|^2
    +\io \big|W''(\chi)\nabla\chi\big|^2
   \le \big(1+\|\nabla \chi\|^2\big)
     Q\big(\|\chi\|_{L^\infty(\Omega)}^2\big)
   \le Q\big(\|\chi\|_{\frac{3+\epsilon}2}^2\big). 
\end{equation}
Next, let us compute \eqref{conto21} plus
$c_1\times$\eqref{energy}. Using also
\eqref{conto22}, we arrive at
\beeq{conto23}
  \ddt\big[\|\chi\|_{\frac{3+\epsilon}2}^2
   +c_1\calE\big]
  +\|\chi_t\|_{\frac{1+\epsilon}2}^2
  \le Q\big(\|\chi\|_{\frac{3+\epsilon}2}^2\big).
\end{equation}
Thus, noting as $\Psi$ the quantity in square brackets
on the \lhs, and using the comparison principle for
ODE's, it follows that there exists a time $T_0>0$,
depending on $\Dzero$ in a monotonically decreasing
way, such that
\beeq{st21}
  \|\Psi\|_{L^\infty(0,T_0)}\le Q(\Dzero),
\end{equation}
whence, integrating 
\eqref{conto23} in time over $(0,T_0)$, 
and recalling \eqref{energycons},
\beeq{st22}
  \|\chi\|_{L^\infty\big(0,T_0;\calV_{\frac{3+\epsilon}2}\big)}
   +\|\chi_t\|_{L^2\big(0,T_0;\calV_{\frac{1+\epsilon}2}\big)}
   \le Q(\Dzero).
\end{equation}
In particular, by the continuous embedding 
$\calV_{\frac{1+\epsilon}2}\subset L^3(\Omega)$, 
we have
\beeq{st23}
  \|\chi_t\|_{L^2(0,T_0;L^3(\Omega))}
   \le Q(\Dzero).
\end{equation}

\vspace{2mm}

\noindent%
{\bf Third estimate.}~~%
Let us note that, since $\teta$ solves
\eqref{calore},
it is $\teta>0$ a.e.~in $\Omega\times (0,\infty)$ 
and $\teta=1$ a.e.~on $\Gamma\times (0,\infty)$.
Thus, we can test \eqref{calore} 
by $\teta^{p-1}-1$ ($p$ given by \eqref{epsip}),
which (at least in an approximation)
lies in $\calV_0$ for 
a.e.~$t\in(0,\infty)$. We then get
\beeq{conto31}
  \ddt\io\Big(\frac1{p}\teta^{p}
   -\teta\Big)
   +\frac{4(p-1)}{(p-2)^2}\big\|\nabla \teta^{\frac{p-2}2}\big\|^2
   \le -\io \big(\teta^{p-1}-1)\chi_t
\end{equation}
and we estimate the \rhs\ as follows:
\bealo 
  -\io \big(\teta^{p-1}-1)\chi_t
   & \le \|\chi_t\|_{L^3(\Omega)} 
    \big\|\teta^{\frac{p-2}2}\big\|_{L^6(\Omega)}
    \big\|\teta^{\frac{p}2}\big\|
   +c\big(1+\|\chi_t\|^2\big)\\[1mm]
 \no
   & \le \frac\sigma{p} \big\|\teta^{\frac{p-2}2}\big\|_{L^6(\Omega)}^2
   +c_\sigma p\|\chi_t\|_{L^3(\Omega)}^2
    \big\|\teta^{\frac{p}2}\big\|^2
   +c\big(1+\|\chi_t\|^2\big)\\[1mm]
 \label{conto32}  
  & \le \frac\sigma{p} \big\|\nabla\teta^{\frac{p-2}2}\big\|^2
   +\frac\sigma{p}
   +c_\sigma p\|\chi_t\|_{L^3(\Omega)}^2
    \|\teta\|_{L^{p}(\Omega)}^{p}
   +c\big(1+\|\chi_t\|_{L^3(\Omega)}^2\big),
\end{align}
where $\sigma>0$ denotes a ``small'' 
constant, independent of $p$, to be chosen
at the end, and correspondingly
$c_\sigma>0$ depends on the same quantities
as the generic $c$ and, additionally,
on the final choice of $\sigma$. In fact, 
passing from row to row, we allow 
$\sigma$ to ``incorporate'' embedding constants.
We used here the continuous embedding 
$V\subset L^6(\Omega)$ and the Young
and Poincar\'e inequalities. Although $p$ is 
a fixed value (cf.~\eqref{epsip}),
here and below we emphasize
the dependence on $p$ of the estimates, since
they will be readily repeated with different
exponents.

Adding 2$\times$\eqref{energy1}
(where the term on the \rhs\ is split
via Young's inequality) to~\eqref{conto31},
multiplying the result by $p$,
and taking $\sigma$ small enough, 
we then obtain
\beeq{conto33a}
  \ddt\io\big[ \teta^{p}
   +p(\teta-2\log\teta) \big]
   +\kappa\big\|\nabla\teta^{\frac{p-2}2}\big\|^2
  \le cp + p^2 c_2 \|\chi_t\|_{L^3(\Omega)}^2
    \big(1+\|\teta\|_{L^{p}(\Omega)}^p\big)
\end{equation}
for some $c_2>0$.
Now, let us set
\beeq{conto63preposto}
  m:=c_2\|\chi_t\|_{L^3(\Omega)}^2,
   \quext{so that, by \eqref{st23}, }\,
   \|m\|_{L^1(0,T_0)}\le Q(\Dzero).
\end{equation}
Defining $\calY$ as $1$ plus the integral on the \lhs\ of
\eqref{conto33a}, we then have
\beeq{conto33c}
  \ddt\calY
   +\kappa\big\|\nabla\teta^{\frac{p-2}2}\big\|^2
  \le cp + p^2m \big(1+\|\teta\|_{L^p(\Omega)}^p\big)
  \le cp + p^2m \calY,
\end{equation}
whence, recalling \eqref{st23} and \eqref{defiX} and
using Gronwall's Lemma,
\beeq{st31}
  \|\teta\|_{L^\infty\big(0,T_0;L^{p}(\Omega)\big)}
   \le Q(\Dzero).
\end{equation}

\vspace{2mm}

\noindent%
{\bf Fourth estimate.}~~%
We test \eqref{calore} by $2t\teta_t/\teta^2$; next, we differentiate
\eqref{phase} in time and test the result by $2t\chi_t$. Taking
the sum and noting that two terms cancel, 
we get
\beeq{conto41}
  \ddt\Big(t\|\chi_t\|^2
   +t\Big\|\nabla\frac1\teta\Big\|^2\Big)
   +2t\io\frac{\teta_t^2}{\teta^2}
   +2t\|\chi_t\|_V^2
   \le (1+2\lambda t)\|\chi_t\|^2
    +\Big\|\nabla\frac1\teta\Big\|^2.
\end{equation}
Then, integrating over $(0,T_0)$ and using 
\eqref{energycons}, we obtain
\beeq{st41}
  \big\|\chi_t(T_0)\big\|^2
   +\Big\|\nabla\frac1{\teta(T_0)}\Big\|^2
   \le Q(\Dzero)\Big(1+\frac1{T_0}\Big)\le Q(\Dzero);
\end{equation}
actually, $(T_0)^{-1}$ depends increasingly on $\Dzero$.

\vspace{2mm}

\noindent%
{\bf Fifth estimate.}~~%
Up to now, we got uniform bounds in the (small) time
interval $[0,T_0]$. Our aim is now to get uniform 
estimates on $[T_0,\infty)$. First, we essentially
repeat the previous estimate, but without the weight $t$.
This gives, of course,
\beeq{conto51}
  \ddt\Big(\|\chi_t\|^2
   +\Big\|\nabla\frac1\teta\Big\|^2\Big)
   +2\io\frac{\teta_t^2}{\teta^2}
   +2\|\chi_t\|_V^2
   \le 2\lambda\|\chi_t\|^2,
\end{equation}
whence, integrating over $(T_0,t)$ for arbitrary $t\ge T_0$,
\bealo
  & \big\|\chi_t(t)\big\|^2
   +\Big\|\nabla\frac1{\teta(t)}\Big\|^2
   +2\int_{T_0}^t\Big(
     \big\|(\log\teta)_t(s)\big\|^2
     +\big\|\chi_t(s)\big\|_V^2\Big)\dis\\
 \label{st51}
  & \mbox{}~~~~~
   \le\big\|\chi_t(T_0)\big\|^2
    +\Big\|\nabla\frac1{\teta(T_0)}\Big\|^2   
    +2\lambda\int_{T_0}^t\|\chi_t(s)\|^2\,\dis
   \le Q(\Dzero),
\end{align}
where we used \eqref{st41} and \eqref{energycons} to
control the terms on the \rhs.

\vspace{2mm}

\noindent%
{\bf Sixth estimate.}~~%
We repeat the Third estimate restarting from $T_0$.
Let us notice that, by \eqref{conto63preposto},
\eqref{st51} and the continuous embedding
$V\subset L^6(\Omega)$, 
\beeq{conto63}
  m=c_2\|\chi_t\|_{L^3(\Omega)}^2
   \quext{satisfies now }\,\,
   M:=\|m\|_{L^1(0,\infty)}\le Q(\Dzero).
\end{equation}
Thus, by the continuous embedding
$V\subset L^6(\Omega)$, \eqref{conto33c} takes the form
\beeq{conto33cnew}
  \ddt \calY
   +\kappa\|\teta\|_{L^{3p-6}(\Omega)}^{p-2}
  \le cp + p^2m \calY.
\end{equation}
Let us now notice that, thanks to $p>3$, we 
can write
\bealo
  \|\teta\|_{L^p(\Omega)}^{p-2}
   & = \Big(\io \teta^p\Big)^{\frac{p-2}p}
     \le \Big(
        \big\|\teta^p\big\|_{L^{\frac{3p-6}p}(\Omega)}
        \big\|1\|_{L^{\frac{3p-6}{2p-6}}(\Omega)}
          \Big)^{\frac{p-2}p}\\
 \label{espop}
   & =  \|\teta\|_{L^{3p-6}(\Omega)}^{p-2}
        |\Omega|^{\frac{2p-6}{3p}}
     \le c_\Omega \|\teta\|_{L^{3p-6}(\Omega)}^{p-2}.
\end{align}
Moreover, we have
\beeq{nuovap}
  \|\teta\|_{L^p(\Omega)}^{p-2}
   = \Big( \calY - 1 - p \io (\teta-2\log\teta) \Big)
          ^{\frac{p-2}p}
   \ge c \calY^{\frac{p-2}p} - p Q(\Dzero).
\end{equation}
Thus, 
\eqref{conto33cnew} gives
\beeq{conto63.1}
  \ddt \calY
   +\kappa\calY^{\frac{p-2}p}
    \le p^2 m \calY + p Q(\Dzero),
\end{equation}
so that, for $\calH:=\log\calY$, 
\beeq{conto63.2}
  \ddt\calH
   + e^{-\frac{2\calH}p}
    \big[ \kappa - p Q(\Dzero) e^{-\frac{(p-2)\calH}p} \big]
   \le p^2 m.
\end{equation}
Noting as $\Sigma$ the quantity in square brackets, 
an easy computation shows that
\beeq{conto63.3}
  \Sigma \ge 0 \Leftrightarrow 
   \calH \ge \frac{p}{p-2} \log
         \Big(\frac{p}\kappa Q(\Dzero)\Big)
           =: \zeta.
\end{equation}
Consequently, it is not difficult to obtain
\beeq{conto63.4}
  \| \calH \|_{L^\infty(T_0,\infty)}
    \le \max\big\{\zeta,\calH(T_0)\big\}
     +p^2\int_0^\infty m(s)\,\dis,
\end{equation}
so that, being by \eqref{conto63.3} 
and \eqref{st31},
\beeq{pezzimmezzo}
  \exp(\zeta) \le p Q(\Dzero), \qquad 
   \exp(\calH(T_0)) 
   = \calY(T_0) \le p Q(\Dzero),
\end{equation}
and using \eqref{conto63}, we readily get
\bealo
   \|\teta\|_{L^\infty(T_0,\infty;L^p(\Omega))}^p  
   & \le \| \calY \|_{L^\infty(T_0,\infty)}
    \le \exp\big( \max \{ \zeta, \calH(T_0) \} \big)
     \exp\big( p^2 M \big)\\
 \label{conto63.5}
   & \le p Q(\Dzero) \exp \big(p^2 M \big),
\end{align}
whence, clearly,
%
\beeq{stimatetap}
  \|\teta\|_{L^\infty(T_0,\infty;L^p(\Omega))} 
   \le Q(\Dzero).  
\end{equation}
Suitable time integrations of \eqref{conto33cnew}
permit us collect what we have proved so far in a
\bele\label{lemmapasso1}
 Under the assumptions of\/ {\rm Theorem~\ref{teounif}},
 there exist a time $T_0>0$ and a quantity $Q_0>0$, both
 depending on $\Dzero$, such that, for all $t\ge T\ge T_0$,
 \beal{primaiter} 
   \|\teta(t)\|_{L^p(\Omega)}^p
    +\int_t^{t+1}\|\teta(s)\|_{L^{3p-6}(\Omega)}^{p-2}\,\dis
    & \le Q_0,\\
  \label{primaiter2}
   \int_{T}^t\|\teta(s)\|_{L^{3p-6}(\Omega)}^{p-2}\,\dis
    & \le Q_0+Q_0(t-T).
 \end{align}
\enle
Our next aim is to extend \eqref{primaiter} and
\eqref{primaiter2} to any finite exponent.
Namely, we have
\bele\label{lemmapasso2}
 Under the assumptions of\/ {\rm Theorem~\ref{teounif}},
 for all $q\in[p,\infty)$ there exist a time $T_q>0$ 
 and a quantity $Q_q>0$, both depending on $\Dzero$ and $q$,
 such that, for all $t\ge T\ge T_q$,
 \beal{secondaiter} 
   \|\teta(t)\|_{L^q(\Omega)}^q
    +\int_t^{t+1}\|\teta(s)\|_{L^{3q-6}(\Omega)}^{q-2}\,\dis
    & \le Q_q,\\
  \label{secondaiter2}
   \int_{T}^t\|\teta(s)\|_{L^{3q-6}(\Omega)}^{q-2}\,\dis
    & \le Q_q+Q_q(t-T).
 \end{align}
\enle
\begin{proof}
 It suffices to iterate finitely many times the procedure
 in the Sixth Estimate. Namely, setting $p_0:=p$, we observe that,
 by \eqref{primaiter} and interpolation, there follows
 \beeq{conto71}
  \sup_{t\in[T_0,\infty)}
   \| \teta \|_{L^{p_i}(t,t+1;L^{p_{i}}(\Omega))}
    \le Q(\Dzero),
 \end{equation}
 where we have set (for $i=1$, at least in the meanwhile)
 \beeq{defpipiu1}
   p_{i}:= \frac{5p_{i-1}-6}3.
 \end{equation}
 Note that $p_1>p_0$ since $p_0>3$.
 Then, we repeat the argument leading to 
 \eqref{conto33cnew}, but with $p_1$ in place 
 of $p$. This gives (for $i=1$ and 
 with obvious meaning of $\calY_i$)
 \beeq{conto72}
   \ddt \calY_i
    +c\|\teta\|_{L^{3p_i-6}(\Omega)}^{p_i-2}
   \le cp_i + p_i^2m \calY_i.
 \end{equation}
 Noting that both $m$ and $\calY_i$ are summable on 
 time intervals of finite length
 thanks to \eqref{conto63}
 and, respectively, \eqref{conto71},
 we can use the {\sl uniform}\/ Gronwall Lemma 
 (cf., e.g., \cite[Lemma III.1.1]{Te}), 
 that gives
 \beeq{conto73}
   \|\teta(t)\|_{L^{p_i}(\Omega)}^{p_i}
    \le \calY_i(t)
    \le Q_i(\Dzero)
     \qquad\perogni t\ge T_i:=T_{i-1}+1,
 \end{equation}
 with obvious meaning of $Q_i(\Dzero)$. 
 Thus, suitable integrations in time 
 of \eqref{conto72} give the analogue of 
 \eqref{primaiter} and \eqref{primaiter2}
 with $p_1$ in place of $q$.
 To get \eqref{secondaiter} and 
 \eqref{secondaiter2}, it 
 then suffices to proceed by iteration
 on $i$ until $p_i$ is larger than $q$. 
 Notice that since a finite number of steps is 
 sufficient, we do not have to take care of 
 the dependence on $i$ of the quantities
 $Q_i$ and $T_i$ (both, in fact, would explode
 if infinite iterations were needed).
 The proof of the Lemma is concluded.
\end{proof}
A similar property holds also for the inverse temperature:
\bele\label{lemmapasso2.5}
 Setting $u:=\teta^{-1}$,
 under the assumptions of\/ {\rm Theorem~\ref{teounif}},
 for all $q\in[1,\infty)$ there exist a time $T_q>0$
 and a quantity $Q_q>0$, both depending on $\Dzero$ and $q$
 (and possibly larger from those in the previous\/
 {\rm Lemma}) such that, for all $t\ge T\ge T_q$,
 \beal{secondaemezza} 
   \|u(t)\|_{L^q(\Omega)}^q
    +\int_t^{t+1}\|u(s)\|_{L^{3q+6}(\Omega)}^{q+2}\,\dis
    & \le Q_q,\\
  \label{secondaemezza2}
   \int_{T}^t\|u(s)\|_{L^{3q+6}(\Omega)}^{q+2}\,\dis
    & \le Q_q+Q_q(t-T).
 \end{align}
\enle
\begin{proof}
 Note that we already know the bound of the first term 
 in \eqref{secondaemezza} for $q=6$ thanks
 to \eqref{st51} and the continuous embedding 
 $V\subset L^6(\Omega)$.
 Then, we proceed essentially as in the 
 Third estimate, i.e., for a generic $q\ge 6$,
 we multiply \eqref{calore}
 by $1-u^{q+1}$. In place of \eqref{conto31},
 we get
 \beeq{conto31u}
   \ddt\io\Big(\frac1{q}u^{q}
    +\teta\Big)
    +\frac{4(q+1)}{(q+2)^2}\big\|\nabla u^{\frac{q+2}2}\big\|^2
    \le - \io \big(1-u^{q+1})\chi_t,
 \end{equation}
 so that, estimating the \rhs\ as in 
 \eqref{conto32}, we infer
 \beeq{conto33cu}
   \ddt\calZ_q 
    +c\|u\|_{L^{3q+6}(\Omega)}^{q+2}
   \le cq + q^2m \calZ_q,
    \qquext{where }\,
    \calZ_q:=\io\big(u^q
     +q\teta\big)
 \end{equation}
 and $m$ is as in \eqref{conto63preposto} (possibly
 for a different value of $c_2$).
 At this point, noticing that the exponents are 
 even better than in \eqref{conto33cnew}, the proof
  can be completed by mimicking the arguments
 in the Sixth estimate and in the proof
 of Lemma~\ref{lemmapasso2}.
\end{proof}
\bele\label{lemmapasso3}
 Under the assumptions of\/ {\rm Theorem~\ref{teounif}},
 there exist a time $T_*$ and a quantity $Q_*$,
 both depending on $\Dzero$, such that, 
 for all $t\ge T_*$,
 \beeq{regochit}
   \|\chi_t(t)\|^2_{L^{24/5}(\Omega)}\le Q_*.
 \end{equation}
\enle
\begin{proof}
 The exponent $24/5$ in \eqref{regochit} is chosen 
 just for later convenience. In fact, \eqref{regochit} 
 can be proved for any exponent strictly smaller than $6$.
 Differentiating in time \eqref{phase}, we have
 \beeq{phaset}
   \chi_{tt}+B\chi_t
    =\Phi:=-W''(\chi)\chi_t+\frac{\teta_t}{\teta^2}
 \end{equation}
 and we claim that, for any $\nu\in(0,1)$,
 we can choose $T_\nu>0$ and $Q_\nu>0$,
 both depending only on $\Dzero$ and $\nu$,
 such that 
 \beeq{stimaPhi}
   \int_{T_\nu}^\infty
    \|\Phi(s)\|^2_{H^{-\nu}(\Omega)}\,\dis
     \le Q(\Dzero).
 \end{equation}
 Actually, recalling \eqref{st51} and applying
 standard regularity results to \eqref{phase}
 (seen here as a time-dependent family of 
 elliptic equations), it follows that 
 \beeq{stimaPhi1}
   \|\chi(t)\|_{H^2(\Omega)}\le Q(\Dzero)
    \quext{for all }\,t\ge T_0.
 \end{equation}
 Thus, by \eqref{W1}, the continuous embedding
 $H^2(\Omega)\subset L^\infty(\Omega)$, and
 \eqref{energycons},
 \beeq{stimaPhi2}
   \int_{T_0}^\infty \|W''(\chi(s))\chi_t(s)\|^2\,\dis
    \le Q(\Dzero).
 \end{equation}
 Analogously, using the first integral bound
 in \eqref{st51}, the bound of the 
 first term in \eqref{secondaemezza}
 with $q$ sufficiently large (depending on $\nu$),
 and elementary interpolation, it is not 
 difficult to get, for some $T_\nu'>0$ and $Q_\nu>0$,
 \beeq{stimaPhi3}
   \int_{T_\nu'}^\infty \Big\|\frac{\teta_t(s)}{\teta^2(s)}
           \Big\|^2_{L^{\frac6{3+2\nu}}(\Omega)}\,\dis
    \le Q_\nu(\Dzero).
 \end{equation}
 Thus, thanks to \eqref{stimaPhi2}, \eqref{stimaPhi3}
 and the continuous embedding
 $L^{\frac6{3+2\nu}}(\Omega)\subset H^{-\nu}(\Omega)$,
 we see that \eqref{stimaPhi} holds for any $T_\nu\ge T_\nu'$.
 Now, let us observe that, by the bound of the second integral 
 term in \eqref{st51}, $T_\nu\in [T_\nu',T_\nu'+1]$ 
 can be chosen such that
 \beeq{regoinizepsi}
   \|\chi_t(T_\nu)\|_{H^{1-\nu}(\Omega)}^2
    \le c\|\chi_t(T_\nu)\|_V^2
    \le Q(\Dzero).
 \end{equation}
 %
 Then, applying the standard
 linear parabolic Hilbert  theory 
 to the equation \eqref{phaset} on the
 time interval $[T_\nu,\infty)$ 
 and with the initial condition $\chi_t(T_\nu)$,
 and using \eqref{stimaPhi}, we have 
 (possibly for a different value of $Q_\nu$)
 \beeq{quasifinale}
   \|\chi_t\|_{C^0([T_\nu,\infty);H^{1-\nu}(\Omega))}
    +\|\chi_t\|_{L^2(T_\nu,\infty;H^{2-\nu}(\Omega))}
   \le Q_\nu(\Dzero),
 \end{equation}
 whence the assert follows from the continuous 
 embedding $H^{1-\nu}(\Omega)\subset L^{24/5}(\Omega)$, 
 which holds for $\nu$ small enough.
\end{proof}

\noindent%
{\bf End of proof of Theorem~\ref{teounif}.}~~
We use a modified Alikakos-Moser \cite{Al} iteration scheme 
similar to that in \cite{La1}, but suitably adapted
in order to obtain time regularization effects. 
Similar procedures have been proved to be effective
in other recent papers, cf.~\cite{EZ,SS}. 

As a first step, we come back to
\eqref{conto31}, where the exponent $p$ is
substituted by a number $q_i$ to be chosen
later. Since we need infinitely many iterations, 
now the \rhs\ has to be estimated more 
carefully. Namely, we have
\bealo 
  -\io \big(\teta^{q_i-1}-1)\chi_t
   & \le \|\chi_t\|_{L^{24/5}(\Omega)} 
    \big\|\teta^{\frac{q_i-2}2}\big\|_{L^6(\Omega)}
    \big\|\teta^{\frac{q_i}2}\big\|_{L^{8/5}(\Omega)}
   +c\big(1+\|\chi_t\|^2\big)\\[1mm]
 \label{conto81}  
  & \le \frac\sigma{q_i} \big\|\nabla\teta^{\frac{q_i-2}2}\big\|^2
   +\frac\sigma{q_i}
   +c_\sigma q_i Q_* \|\teta\|_{L^{4q_i/5}(\Omega)}^{q_i}
   +C_0,
\end{align}
where $Q_*$ is exactly the same quantity as in 
\eqref{regochit} and the constant $C_0$
depends on $\Dzero$ and is independent of $q_i$.
Thus, possibly modifying $C_0$, 
in place of \eqref{conto33cnew} we get
\beeq{conto82}
  \ddt \calY_i
   +\kappa\|\teta\|_{L^{3q_i-6}(\Omega)}^{q_i-2}
  \le c q_i^2 Q_* \|\teta\|_{L^{4q_i/5}(\Omega)}^{q_i}
   + C_0 q_i,
\end{equation}
where it is worth noting that 
\beeq{defiYi}
  \|\teta\|_{L^{q_i}(\Omega)}^{q_i}
  \le \calY_i:=1+\io\big[\teta^{q_i}
   +q_i(\teta-2\log\teta)\big]
     \le\|\teta\|_{L^{q_i}(\Omega)}^{q_i}
     + C_1 q_i,
\end{equation}
where $C_1$ is another quantity depending only on
$\Dzero$. Moreover, \eqref{conto82} gives
\beeq{conto83}
  \ddt \calY_i
   +\kappa\|\teta\|_{L^{3q_i-6}(\Omega)}^{q_i-2}
  \le c q_i^2 Q_* \|\teta\|_{L^{\frac{4q_i}5}(\Omega)}^{\frac{q_i}5+2}
   \|\teta\|_{L^{\frac{12q_i}5-6}(\Omega)}^{\frac{4q_i}5-2}
   + C_0 q_i,
\end{equation}
provided that $12q_i/5-6\ge 4q_i/5$, which is 
true for all $i\in \NN$
if $q_0$ is large enough and we set
\beeq{sceltaqi}
  q_i=\frac54q_{i-1}, \quad \perogni i\ge 1.
\end{equation}
The choice \eqref{sceltaqi} permits to rewrite
\eqref{conto83} in the form
\beeq{conto84}
  \ddt \calY_i
   +\kappa\|\teta\|_{L^{3q_i-6}(\Omega)}^{q_i-2}
  \le c_3 q_i^2 Q_* \calY_{i-1}^{\frac14+\frac5{2q_i}}
   \|\teta\|_{L^{3q_{i-1}-6}(\Omega)}^{q_{i-1}-2}
   + C_0 q_i,
\end{equation}
for some $c_3>0$. Let us now define 
\beeq{defitauinfty}
  \tau_i:=i^{-2}, \qquext{so that }\,\,
   \tau_\infty:=\sum_{i=1}^\infty \tau_i
    =\sum_{i=1}^\infty i^{-2}<\infty.
\end{equation}
Moreover, let us take $T_*$ as in Lemma~\ref{lemmapasso3}
and set $T_\infty:=T_*+\tau_\infty$. We 
now aim to show that 
for all fixed $t\ge T_\infty$ the bounds 
\eqref{rego1}--\eqref{rego2} are satisfied. 
More precisely, we will limit ourselves
to prove the $L^\infty$-bound 
of $\teta$ in \eqref{rego1}. 
Indeed, as noted in the proof of 
Lemma~\ref{lemmapasso2.5},
the argument to prove the $L^\infty$-bound
of $u$ is similar and even simpler;
moreover, the $V$-bounds are consequence
of the $L^\infty$-bounds and of \eqref{st51};
finally, \eqref{rego2} is already known
from~\eqref{stimaPhi1}.

Thus, we assume that $t\ge T_\infty$ is fixed and
set $S_0:=t-\tau_\infty$, so that $S_0\ge T_*$.
We shall now work on the interval 
$[t-\tau_\infty,t-\tau_\infty]=[S_0,S_0+2\tau_\infty]$.
Using Lemma~\ref{lemmapasso2}, we can also 
assume $q_0$ as large as we want, so that
there exists a quantity $R_0$,
depending only on $\Dzero$ and on the choice of 
$\tau_\infty$, such that the bound
\beeq{firstiter}
  \|\calY_0\|_{L^\infty(S_0,s)} 
   +\kappa\int_{S_0}^s \|\teta\|_{L^{3q_{0}-6}(\Omega)}^{q_{0}-2}
    \le R_0,
\end{equation}
where $\calY_0$ is defined as in \eqref{defiYi} with $i=0$,
holds uniformly w.r.t.~$s\in [S_0,S_0+2\tau_\infty]$.
Thus, taking $i=1$, integrating \eqref{conto84} 
in time over $(S_1,s)$,
where $S_1\in[S_0,S_0+\tau_1]=[S_0,S_0+1]$
will be chosen later and $s$ is a generic
point in $[S_1,S_0+2\tau_\infty]$
so that $s-S_1\le 2\tau_\infty$, we have
\beeq{conto85}
  \calY_1(s) 
   +\kappa\int_{S_1}^s\|\teta\|_{L^{3q_1-6}(\Omega)}^{q_1-2}
  \le \calY_1(S_1)
   + \frac{c_3}\kappa \Big(\frac54\Big)^2 q_0^2 Q_* 
         R_0^{\frac54\big(\frac{5q_0+8}{5q_0}\big)} 
   + \frac54 C_0 q_0(2\tau_\infty)
\end{equation}
and we notice that also the first term on the \rhs\
can be estimated. Indeed, by the latter of 
\eqref{firstiter}, $S_1\in[S_0,S_0+\tau_1]$ 
can be chosen such that 
\beeq{stimainiz}
  \|\teta(S_1)\|_{L^{q_1}(\Omega)}^{q_0-2}
    \le \frac1{\tau_1}\int_{S_0}^{S_0+\tau_1}
       \|\teta\|_{L^{q_1}(\Omega)}^{q_0-2}
    \le c\frac{\kappa}{\tau_1}\int_{S_0}^{S_0+\tau_1}
       \|\teta\|_{L^{3q_0-6}(\Omega)}^{q_0-2}
    \le c_4\frac {R_0}{\tau_1}=c_4R_0,
\end{equation}
for some $c_4>0$. We used that
$\tau_1=1$. Recalling \eqref{defiYi}, 
as a consequence we obtain
\beeq{stimainiz2}
  \calY_1(S_1)
   \le \|\teta(S_1)\|_{L^{q_1}(\Omega)}^{q_1}+q_1C_1
   \le \Big(c_4\frac{R_0}{\tau_1}\Big)^{\frac54\frac{q_0}{q_0-2}}
      +\frac54q_0 C_1
   = (c_4 R_0)^{\frac54\frac{q_0}{q_0-2}}+\frac54q_0 C_1.
\end{equation}
Then, setting for $i\ge0$
\beeq{defieta1}
  \eta_i:=\frac{q_i}{q_i-2}
    \ge \frac{5q_i+8}{5q_i}
\end{equation}
and collecting \eqref{conto85}--\eqref{stimainiz2},
we obtain that for all $s\in [S_1,S_0+2\tau_\infty]$,
\beeq{conto86}
  \calY_1(s) 
   +\kappa\int_{S_1}^s\|\teta\|_{L^{3q_1-6}(\Omega)}^{q_1-2}
  \le R_0^{\frac54\eta_0}\Big(
    c_4^{\frac54\eta_0}
    + \frac{25}{16}\frac{c_3}{\kappa}q_0^2 Q_* \Big)
    + \frac54q_0\big(2\tau_\infty C_0+C_1)
  =:R_1.
\end{equation}
At this point we can proceed by iteration
and observe that, as the procedure is repeated,
the main modification comes from a term
$t_i^{-1}=i^2$ additionally appearing on
the \rhs\ of the $i$-analogue of~\eqref{stimainiz}.
Suitably modifying the procedure,
\eqref{conto86} takes the new form
\beeq{iiter}
  \|\calY_i\|_{L^\infty(S_i,s)} 
   + \kappa \int_{S_i}^s \|\teta\|_{L^{3q_{i}-6}(\Omega)}^{q_{i}-2}
  \le R_i, \qquad \perogni s\in [S_i,S_0+2\tau_\infty],
\end{equation}
where we point out that $\kappa$, which comes
from \eqref{conto84} and, in fact,
from \eqref{conto33cnew}, is independent of $i$.
Moreover, $S_i$ is a suitable point 
in $[S_{i-1},S_{i-1}+\tau_i]$ 
and $R_i$ is given by
\bealo
  R_i & = R_{i-1}^{\frac54\eta_{i-1}}\Big(
   (c_4i^2)^{\frac54\eta_{i-1}}
    +\Big(\frac{25}{16}\Big)^i\frac{c_3}{\kappa}q_0^2 Q_*\Big)
    + \Big(\frac54\Big)^iq_0\big(2\tau_\infty C_0+C_1)\\
 \label{Riipiu1}
   & \le \Big((c_4i^2)^{\frac54\eta_{i-1}}
   +\Big(\frac{25}{16}\Big)^i K\Big)
    R_{i-1}^{\frac54\eta_{i-1}}
    =: A_i R_{i-1}^{\frac54\eta_{i-1}},
\end{align}
for some $K>0$ depending on $\Dzero$ and 
the choice of $\tau_\infty$, and independent of $i$.
Consequently, $R_i$ is estimated in terms of $R_0$ by
\beeq{stimaRiR0}
  R_i \le R_0^{(\frac54)^i\prod_{k=0}^{i-1}\eta_k}
   A_i A_{i-1}^{\frac54\eta_{i-1}}
    A_{i-2}^{(\frac54)^2\eta_{i-1}\eta_{i-2}}
    \dots
   = R_0^{(\frac54)^i\prod_{k=0}^{i-1}\eta_k}
   \prod_{j=1}^i
   A_j^{(\frac54)^{i-j}\prod_{k=j}^{i-1}\eta_k},
\end{equation}
where it is intended that the latter productory
is $1$ in the case $j=i$. Passing to the 
logarithm and observing that
\beeq{Rilimi}
  \prod_{k=1}^{\infty}\eta_k<\infty,
\end{equation}
it is then easy to verify that
\beeq{Rilimi2}
  \limsup_{i\nearrow\infty} R_i^{\frac1{q_i}}
   = \limsup_{i\nearrow\infty} R_i^{(\frac45)^i\frac1{q_0}}
   \le Q(\Dzero) < \infty.
\end{equation}
Thus, coming back to \eqref{stimaRiR0},
recalling \eqref{defiYi}, and noting that the sequence
$S_i$ converges to a point $S_\infty$ such that 
$S_0\le S_\infty\le t=S_0+\tau_\infty \le S_0+2\tau_\infty$,
we finally infer that
\beeq{Rilimi3}
  \limsup_{i\nearrow\infty} \|\teta(s)\|_{L^{q_i}(\Omega)}
    \le   \limsup_{i\nearrow\infty} R_i^{\frac1{q_i}}
    \le Q(\Dzero) \qquad \perogni s\in[S_\infty,S_0+2\tau_\infty].
\end{equation}
In particular, this holds for $s=t$ and
it is worth remarking once more that the latter
quantity $Q(\Dzero)$ is independent of 
$t\in[T_*+\tau_\infty,\infty]$. Actually,
it depends on time only through the choice of the
sequence $\tau_i$, and not on the choice of $S_0\ge T_*$,
i.e., of $t$. The proof of the first of \eqref{rego1} 
and of the Theorem is complete.\dimbox


\vspace{2mm}

\noindent%
{\bf Proof of Theorem~\ref{teoattra}.}~~%
We start noticing that $\calE$, defined in \eqref{defiE},
is a Liapounov functional for Problem~(P). Namely, the 
following conditions (cf., e.g., \cite[Sec.~5]{Ba1})
hold:\\[2mm]
 (L1)~~%
 $\calE$ is continuous on $\calX$ 
 (recall~\eqref{defiX});\\[1mm]
 (L2)~~%
 $\calE$ is nonincreasing along solution 
 trajectories;\\[1mm]
 (L3)~~%
 if $S(t)w=w$ for some $t>0$ and $w\in \calX$, then
 $w$ belongs to the set $\calEq$ of {\sl equilibrium
 points}\/ of the semigroup (and consequently it identifies
 a stationary solution).\\[2mm]
Indeed, (L1) is obvious since $\calX$ is endowed with
the metric~\eqref{defidX}; (L2) is a simple consequence
of the energy equality \eqref{energy}; finally, (L3)
still follows from \eqref{energy} by noticing
that if $\chi_t=0$ and $\teta^{-1}=1$ then it is 
also $\teta_t=0$ by comparison in \eqref{calore}. It is 
worth remarking that $w=(\teta,\chi)\in\calX$ is a 
stationary point of $S(\cdot)$ if and only if
$\teta\equiv 1$ in $\Omega$ and $\chi$ solves
\beeq{statprob}
  B\chi + W'(\chi) = 0, \qquext{in }\,V'.
\end{equation}
It is well-known that, due to nonconvexity
of $W$, \eqref{statprob} can have infinitely 
many solutions, so that the structure of
$\omega$-limits of solutions to (P) and, a 
fortiori, of attractors, is nontrivial. 
Nevertheless, by maximum principle arguments and 
standard elliptic regularity theorems
(cf., e.g., \cite{ADN}), it is easy to 
prove that the projection of $\calEq$ 
on the second component $\chi$ is bounded 
at least in $W^{2,\zeta}(\Omega)$ for all 
$\zeta\in [1,\infty)$ (actually, using bootstrap
arguments, more could be said depending
on the smoothness of $W$, but we are not
interested in maximal regularity here).

Thus, let $\calB_0$ be the neighbourhood
of $\calEq$ of radius $1$ in the metric
of $\calX$. Then, a simple and direct
contradiction argument (cf.~\cite[Thm.~5.1]{Ba1}) 
shows that $\calB_0$ is pointwise 
absorbing for $S(\cdot)$, i.e.,
given a solution $w$ to (P)
with initial datum in $\calX$, there
exists a time $T_w$ such that 
$w(t)\in \calB_0$ for all $t\ge T_w$.
Thus, being $S(\cdot)$ {\sl asymptotically
compact}\/ (i.e., $S(\cdot)$ eventually
maps $\calX$-bounded sets of initial data 
into relatively compact sets) thanks to
\eqref{rego1}--\eqref{rego2}, we deduce
existence of the global attractor~$\calA$
by means, e.g., of \cite[Thm.~3.3]{Ba1}.

To complete the proof, we have to show that
\eqref{compA} holds. To do this, it suffices
to notice that, as a consequence of the existence
of $\calA$, $S(\cdot)$ admits a $\calX$-bounded
and {\sl uniformly}\/ absorbing set $\calB_1$.
Namely, for every $\calX$-bounded set $B$ there
exists $T_B$ such that $S(t)B\subset \calB_1$ for
all $t\ge T_B$. We then notice that,
by Theorem~\ref{teounif}, 
for sufficiently large $t$, $S(t)$ maps
$\calB_1$ into a set $\calB_2$
which is bounded in the same sense as \eqref{compA}.
Thus, $\calB_2$ is absorbing because
$\calB_1$ is absorbing, and, consequently,
the bound \eqref{compA} holds also 
for the attractor $\calA$ which is the 
$\omega$-limit of $\calB_2$. The proof
is concluded.\dimbox

\vspace{2mm}

\noindent%
{\bf Proof of Theorem~\ref{teoespo}.}~~%
Let us recall the basic uniqueness estimate for
system \eqref{calore}--\eqref{phase}. Let 
$(\teta_i,\chi_i)$, $i=1,2$, be a couple 
of solutions to~(P) and set 
$(\teta,\chi):=(\teta_1,\chi_1)-(\teta_2,\chi_2)$.
Set also $e_i:=\teta_i+\chi_i$, $i=1,2$, 
and $e:=e_1-e_2$ (the new variable has the physical
meaning of {\sl enthalpy}).
Write the differences of \eqref{calore} and 
\eqref{phase} for $i=1,2$ and test them, respectively, 
by $A^{-1}e$ and by $\chi$. Taking the sum,
noting that two terms cancel, and using 
\eqref{W2}, we then obtain
\beeq{stuniq}
  \ddt\big(\|e\|_{H^{-1}(\Omega)}^2
   +\|\chi\|^2\big)
   +2\io \Big(-\frac1{\teta_1}+\frac1{\teta_2}\Big)\teta
   +2\|\chi\|_V^2
   \le 2\lambda \|\chi\|_H^2.
\end{equation}
%

Now, let us restrict ourselves to consider only
initial data lying in a suitable absorbing set.
Namely, we take the absorbing set $\calB_2$ defined
above and set
\beeq{defiB3}
  \calB_3:=\overline{\{\cup_{t\ge T_2}S(t)\calB_2\}},
\end{equation}
where $T_2>0$ is such that $\calB_2$ absorbs itself for $t\ge T_2$
and we have taken what we will call
the {\sl sequential weak star}\/ closure in $\calW$.
Namely, we define the set $\calW$~as
\beeq{defiW}
  \calW:=\big\{(\teta,\chi)\in (V\cap L^\infty(\Omega))\times
    H^2(\Omega):~\teta>0~\text{a.e., and}~ 
    u=\teta^{-1}\in V\cap L^\infty(\Omega)\big\}
\end{equation}
and we intend that a point $w=(\teta,\chi)$ belongs 
to $\calB_3$ iff there
exist two sequences $\{w_n\}=\{(\teta_n,\chi_n)\}\subset \calB_2$
and $\{t_n\}\subset [T_2,\infty)$ such that, as $n\nearrow\infty$,
$S(t_n)w_n=:(\teta_n(t_n),\chi_n(t_n))$ satisfies
\beeq{Wweakstar}
  \teta_n(t_n)\to \teta,~~\teta_n^{-1}(t_n)\to \teta_{-1},
   ~~\text{weakly-$*$ in }\,V\cap L^\infty(\Omega)
    \quext{and }\, \chi_n(t_n)\to \chi,~~
         \text{weakly in }\,H^2(\Omega).
\end{equation}
Of course, this {\sl weak star}\/ convergence of 
$\calW$ is associated to 
a suitable (Hausdorff) topology, which we note
as the ``weak star topology'', or simply 
the ``topology'' of $\calW$.
Instead, when we speak, e.g., of the 
$(H\times V)$-norm of an element
$w=(\teta,\chi)\in\calW$, we will just mean
$(\|\teta\|^2+\|\chi\|_V^2)^{1/2}$ so that
we are neglecting, in fact, the behavior 
of $u=\teta^{-1}$. 
Thus, the generic element of $\calW$ is 
seen just as a couple; however,
the ``weak star convergence'' defined in 
\eqref{Wweakstar} and the related topology
take also the additional variable $u$ 
into account.

Next, it is worth noticing that, 
by construction, $\calB_3$ is 
positively invariant (i.e.~$S(\tau)\calB_3\subset \calB_3$
for all $\tau\ge 0$), sequentially weakly star closed
in~$\calW$, and contained in the 
$\calW$-sequential weak star closure of $\calB_2$,
so that, in particular, there holds
(cf.~\eqref{rego1}--\eqref{rego2})
\beeq{propB3}
  \|\teta\|_{V\cap L^\infty(\Omega)}
   +\|\teta^{-1}\|_{V\cap L^\infty(\Omega)}
   +\|\chi\|_{H^2(\Omega)}
   \le \cc_3 
\end{equation}
for all $(\teta,\chi)\in\calB_3$
and for some constant $\cc_3>0$.
%
%
%
Consequently, we have
\beeq{stuniq2}
  2\io \Big(-\frac1{\teta_1}+\frac1{\teta_2}\Big)\teta
   \ge c_5 \|\teta\|^2
   \qquad\perogni \teta\in \Pi_1(\calB_3).
\end{equation}
where $\Pi_1$ is the projection on the first component
and $c_5$ suitably depends on $\cc_3$.

We now refer to the so-called method of 
$\ell$-trajectories (cf.~\cite{MP,Pr1,Pr2,Pr3}).
To do this, let us take $\ell>0$ and define 
the set $\calUl$ of $\ell$-trajectories 
of~(P) simply as the set of the solutions
whose initial datum lies in $\calB_3$,
restricted to the time interval $(0,\ell)$.
Thus, as before, solutions are seen as 
couples $(\teta,\chi)$ and the behavior
of $u=\teta^{-1}$ is not considered; nevertheless,
since the elements of $\calUl$ take values in 
$\calB_3$ (recall that $\calB_3$ is positively
invariant), they satisfy \eqref{propB3}
uniformly in time.
%
The set $\calUl$ is endowed with the norm
of $L^2(0,\ell;H\times V)$ (the reason for a choice
of such a weak metric is in estimate
\eqref{stuniq}). Let us notice, however, that,
if $\{w_n\}\subset \calUl$ tends 
to some limit $w$ (strongly) in 
$L^2(0,\ell;H\times V)$, then 
also $w$ lies in $\calUl$ (in other words, 
$\calUl$ is complete in the chosen metric).
Indeed, by construction, for 
all $n$ there holds
\beeq{inB3}
  w\zzn:=w_n(0)=
   \lim_{k\to\infty}S(t_k^n)w\zzn^k,
   \qquext{where }\/w\zzn^k\in \calB_2,
\end{equation}
$t_k^n\ge T_2$ for all $k$ and $n$, and the limit is
intended in the topology of
$\calW$ (cf.~\eqref{Wweakstar}). 
In particular, both the above $k$-limit 
and the $n$-limit $w\zzn\to w_0$ (the latter
thanks to \eqref{propB3}) hold in $\calX$, which is 
a metric space. Thus, 
we can extract a diagonal subsequence 
such that
\beeq{inB3-2}
  \lim_{j\nearrow\infty}S(t^{n_j}_{k_j})w_{0,n_j}^{k_j}\to w_0:=w(0),
    \qquext{in }\,\calX
\end{equation}
and, again by uniform validity of \eqref{propB3},
weakly star in 
$\calW$. This means that $w_0\in \calB_3$ and 
$w\in \calUl$, as desired.

Let us now integrate \eqref{stuniq}
over $(\tau,2\ell)$, 
where $\tau$ is a generic point 
in $[0,\ell]$. We then get
\bealo
  & \|e(2\ell)\|_{H^{-1}(\Omega)}^2
   +\|\chi(2\ell)\|^2
   +\int_\tau^{2\ell}
    \big(c_5\|\teta(s)\|^2 
        +2\|\chi(s)\|_V^2\big)\dis\\
 \label{stuniqell}
  & \mbox{}~~~~~
   \le \|e(\tau)\|_{H^{-1}(\Omega)}^2
   +\|\chi(\tau)\|^2
   + 2\lambda \int_\tau^{2\ell}\|\chi(s)\|^2\,\dis.
\end{align}
Integrating the above relation with respect
to $\tau\in(0,\ell)$, we obtain
\bealo
  & \ell\|e(2\ell)\|_{H^{-1}(\Omega)}^2
   +\ell\|\chi(2\ell)\|^2
   +\ell\int_\ell^{2\ell}
    \big(c_5\|\teta(s)\|^2 
        +2\|\chi(s)\|_V^2\big)\dis\\
 \label{stuniqell2}
  & \mbox{}~~~~~
   \le \int_0^\ell c_6\big(\|\teta(\tau)\|_{H^{-1}(\Omega)}^2
   +\|\chi(\tau)\|^2\big)\,\ditau
   + 2\lambda\ell \int_0^{2\ell}\|\chi(s)\|^2\,\dis.
\end{align}
Now, let us use the following
straighforward fact
(see, e.g., \cite[Lemma~3.2]{Pr1}):
\bele\label{daPrazak}
 Let $\calH$ be a Hilbert space and 
 $\calW$ a Banach space such that 
 $\calH$ is compactly embedded into 
 $\calW$. Then, for any $\gamma>0$ 
 there exist a finite-dimensional orthonormal 
 projector $P:\calH\to \calH$ and a positive
 constant $k$, both depending on $\gamma$ and
 such that, for all $z\in \calH$,
 \beeq{compastr}
   \|z\|^2_\calW \le \gamma\|z\|_\calH^2
    +k\|Pz\|_\calH^2.\dimbox
 \end{equation}
\enle
We apply here the Lemma to $z=\chi$ with
$\calH=V$ and $\calW=H$ and to $z=\teta$ with 
$\calH=H$ and $\calW=V_0'=H^{-1}(\Omega)$.
Then, introducing the time shift operator $\calL$, given 
by $\calL:v(\cdot)\mapsto v(\cdot+\ell)$
(where $v$ is a generic function of time),
and dividing \eqref{stuniqell2} by $\ell$,
we obtain
\bealo
  & c_5 \|\calL \teta\|_{L^2(0,\ell;H)}^2 
  +2 \|\calL \chi\|_{L^2(0,\ell;V)}^2
    \le c_6\ell^{-1}\big( \|\teta\|_{L^2(0,\ell;H)}^2 
   +\|\chi\|_{L^2(0,\ell;V)}^2 \big)\\
 \label{newconto56}
  & \mbox{}~~~~~~~~~~
    +2\gamma\lambda \big( \|\calL\chi\|_{L^2(0,\ell;V)}^2 
   +\|\chi\|_{L^2(0,\ell;V)}^2 \big)
   +2k\lambda \big( \|P\calL\chi\|_{L^2(0,\ell;V)}^2 
   +\|P\chi\|_{L^2(0,\ell;V)}^2 \big),
\end{align}
whence, recalling the notation $w:=(\teta,\chi)$ 
and rearranging, 
\bealo
  & \min\big\{c_5,2-2\gamma\lambda\big\}
   \|\calL w\|_{L^2(0,\ell;H\times V)}^2\\
 \label{proviamo}
  & \mbox{}~~~~~
 \le \Big(\frac{c_6}{\ell}+2\gamma\lambda\Big)
   \|w\|_{L^2(0,\ell;H\times V)}^2
   +\cc \big(
  \|Pw\|_{L^2(0,\ell;H\times V)}^2
  +\|P\calL w\|_{L^2(0,\ell;H\times V)}^2
           \big),
\end{align}
where $\cc$ depends on $\gamma, \lambda, \ell$
and all the other constants.
Being not restrictive to assume $c_5\le1$, 
it is clear that we can divide the above
by $c_5$ and choose $\ell$ large enough
and $\gamma$ small enough to obtain (clearly
for a different value of $\cc$)
%
%
%
\beeq{newconto57}
  \|\calL w\|_{L^2(0,\ell;H\times V)}^2
   \le \frac18\|w\|_{L^2(0,\ell;H\times V)}^2
   +\cc \big(\|Pw\|_{L^2(0,\ell;H\times V)}^2
     +\|P\calL w\|_{L^2(0,\ell;H\times V)}^2\big).
\end{equation}
Consequently, the semigroup $S(\cdot)$ 
enjoys the {\sl generalized squeezing
property}\/ introduced in \cite[Def.~3.1]{Pr1}
on the set $\calB_3$. Recalling 
\cite[Lemma~2.2]{Pr2}, 
we then infer that the 
discrete dynamical system on $\calUl$ 
generated by $\calL$ admits 
an exponential attractor
$\calM_{\discr}$.

To conclude, we have to prove that, in 
fact, we can build the exponential attractor
for the original semigroup $S$. Here, however, 
we have to pass from the $H\times V$ to the 
weaker $V_0'\times H$-topology (we recall
that $V_0'=H^{-1}(\Omega)$). Actually,
we can observe that the following
properties hold:\\[2mm]
(M1)~~%
The evaluation map $e:\calUl\to V_0'\times H$ given by
$e:w\mapsto w(\ell)$ is Lipschitz continuous.
%
%
%
%
%
%
%
%
%
%
%
%
To see this, it suffices to multiply \eqref{stuniq}
by $t$ and integrate in $\dit$ between 
$0$ and $\ell$. Notice that, more precisely,
Lipschitz continuity still holds as $\calW_\ell$
is endowed with the weaker
$L^2(0,\ell;V_0'\times H)$-norm;\\[1mm]
(M2)~~%
The map $S(t)$ is uniformly
Lipschitz continuous on $[0,\ell]$ in the
sense that
\beeq{uniLip}
  \|S(t)w_1-S(t)w_2\|_{H^{-1}(\Omega)\times H}
   \le c(\ell) \|w_1-w_2\|_{H^{-1}(\Omega)\times H},
   \qquad\perogni w_1,w_2\in \calB_3~~
    \text{and }\,\perogni t\in [0,\ell].   
\end{equation}
This is easily shown by integrating once 
more \eqref{stuniq} over $(0,t)$ and using the 
Gronwall Lemma.\\[1mm]
(M3)~~%
For each solution $w\in\calUl$ and all 
$0\le s\le t\le \ell$, 
by interpolation there holds
\bealo
  \|w(t)-w(s)\|_{H\times V}^2
  & \le \|w(t)-w(s)\|_{H^1_0(\Omega)\times H^2(\Omega)}
    \|w(t)-w(s)\|_{H^{-1}(\Omega)\times H}\\
 \label{Holderintime}
  & \le C \Big|\int_s^t \|w_t(\tau)\|_{H^{-1}(\Omega)\times H}
  \,\ditau\Big|
   \le C|t-s|^{1/2},
\end{align}
where the constants $C$ depend on the 
``radius'' of $\calB_3$ in $\calW$ (i.e.~on
$\cc_3$, cf.~\eqref{propB3})
and the latter estimate is 
a consequence of the regularity
properties \eqref{regtetaFS}--\eqref{regchiFS}.
Thus, $\ell$-trajectories in $\calUl$ are uniformly 
H\"older continuous in time (notice that this even
holds in the $H\times V$-norm).
Then, properties~(M1)--(M3) allow us to apply, 
e.g., \cite[Thm.~2.6]{MP}, which states that
there exists a set $M$ which is compact and
has finite fractal dimension in 
$V_0'\times H$, is positively
invariant, and exponentially attracts
$\calB_3$. Setting $\calM:=M\cap \calB_3$,
of course $\calM$ satisfies the same properties
of $M$ and, additionally, is bounded in 
$\calW$ (in the sense of \eqref{propB3})
and therefore compact in $\calX$. 
To conclude the proof 
of Theorem~\ref{teoespo}, we have to
show that $\calM$ attracts exponentially fast
any bounded $B\subset \calX$. Actually, this is
true since $\calB_3$ is uniformly absorbing
(so that it exponentially attracts
$B$) and one can use the contractive 
estimate \eqref{stuniq} and 
the {\sl transitivity}\/ property of 
exponential attraction proved in 
\cite[Thm.~5.1]{FGMZ}. Notice 
in particular that the constant
$\kappa$ in \eqref{espoattra}
can be taken independent of $B$ since 
$\lambda$ on the \rhs\ of \eqref{stuniq}
is also independent of $B$
(cf.~\cite[(5.1)]{FGMZ}).
\beos\label{secondaregoespo}
 We can see that exponential attraction 
 still holds in the (stronger)
 $V_0'\times V$-norm.
 This requires to show (M1) and (M2)
 with respect to that topology. To do this, 
 we can write the difference 
 of \eqref{phase} and test it by 
 $\chi_t=\chi_{1,t}-\chi_{2,t}$. Standard
 manipulations then lead to 
 %
 \beeq{stuniq3}
   \ddt\|\chi\|_V^2
    +\|\chi_t\|^2
   \le c_7\big(\|\chi\|^2+\|\teta\|^2\big),
 \end{equation}
 where $\teta=\teta_1-\teta_2$
 and $c_7$ depends on the 
 ``$\calW$-radius'' of $\calB_3$.
 Then, multiplying
 \eqref{stuniq3} by $c_5/2c_7$ 
 ($c_5$ being as in \eqref{stuniq2})
 and summing to \eqref{stuniq},
 we get a contractive estimate in the desired
 topology. On the contrary, 
 at least in the three-dimensional case,
 it seems more difficult to
 obtain an $H$-contraction estimate for $\teta$
 (i.e.~to pass to the $H\times V$-norm).
 Actually one could test the difference of 
 \eqref{calore} by $\teta$. However,
 even knowing the boundedness
 \eqref{propB3}, getting a 
 control of the term involving the Laplacean
 seems out of reach.
\eddos




\vspace{20mm}

\noindent%
{\bf Author's address:}\\[1mm]
Giulio Schimperna\\
Dipartimento di Matematica, Universit\`a degli Studi di Pavia\\
Via Ferrata, 1,~~I-27100 Pavia,~~Italy\\
E-mail:~~{\tt giusch04@unipv.it}

\end{document}